\documentclass[a4paper,10.5pt]{article}
\usepackage{graphicx}
\usepackage{amssymb}
\usepackage{latexsym, bm}
\usepackage{multicol}
\usepackage{indentfirst}
\usepackage{amssymb,amsfonts}
\usepackage{amsmath}
\newtheorem{theorem}{Theorem}
\newtheorem{lemma}{Lemma}

\usepackage{amssymb}
\title{A Folkman Linear Family\thanks{Supported in part by NSFC and CSC.}}
\author{ Qizhong Lin$^a$,\,\,  Yusheng Li$^b$ \vspace*{0.3cm}\\
$^a$Center for Discrete Mathematics, Fuzhou University \\Fuzhou 350108, China\\
$^b$Department of Mathematics, Tongji University  \\Shanghai 200092, China \\
{\small\em Email: linqizhong@fzu.edu.cn, li\_yusheng@tongji.edu.cn}}
\date{}
\begin{document}
\maketitle
\begin{abstract}
For graphs $F$ and $G$, let $F\to (G,G)$ signify that any red/blue edge coloring of $F$ contains a monochromatic $G$. Define Folkman number $f(G;p)$ to be the smallest order of a graph $F$ such that $F\to (G,G)$ and $\omega(F) \le p$.
It is shown that $f(G;p)\le cn$ for graphs $G$ of order $n$ with $\Delta(G)\le \Delta$, where $\Delta\ge 3$, $c=c(\Delta)$ and $p=p(\Delta)$ are positive constants.

\medskip

\end{abstract}

\section{Introduction}

For graphs $F$ and $G$, let $F \to (G,G)$ signify that any red/blue edge coloring of $F$ contains a monochromatic $G$.
The Ramsey number $R(G)$ is the smallest $N$ such that $K_N\to (G,G)$. For most graphs $G$, it is difficult to determine the behavior of $R(G)$, and even more
difficult if the edge-colored graphs are restricted within that of
smaller cliques instead of the complete graphs.

Define a family ${\cal F}(G;p)$ of graphs as
\[
{\cal F}(G;p)=\{F: F\to
(G,G)\hspace{2mm}\mbox{and}\hspace{2mm} \omega(F)\le p\},
\]
where $\omega(G)$ is the clique number of $G$, and define
\[
f(G;p)=\min\{|V(F)|:\; F \in {\cal F}(G;p) \},
\]
which is called the Folkman number. We admit that $f(G;p)=\infty$ if ${\cal F}(G;p)=\emptyset$, and thus
$f(G;p)=\infty$ if $p< \omega(G)$.

The investigation was motivated by a question of Erd\H{o}s and
Hajnal \cite{erd-haj} who asked what was the minimum $p$ such that
${\cal F}(K_3;p)\not=\emptyset$. An important result of Folkman
\cite{folkman} states that ${\cal F}(K_n;p)\not=\emptyset$ for $p\ge n$, which was generalized by Ne\v{s}et\v{r}il and
R\"{o}dl \cite{nes-rodl} as ${\cal F}(G;p)\not=\emptyset$ for $p\ge \omega(G)$. The following property is clear.

\begin{lemma}\label{lem-dec}
The function $f(G;p)$ is decreasing on $p$, and if $p\ge R(G)$, then $f(G;p)=R(G)$.
\end{lemma}

Graham \cite{graham} proved that $f(K_3;5)= 8$ by showing $K_8\setminus C_5\not\to (K_3,K_3)$.
Irving \cite{irving} proved that $f(K_3;4)\le 18$, and it was further improved by Khadzhiivanov
and Nenov \cite{k-n} to $f(K_3;4)\le 16$. Finally, Piwakowski, Radziszowski, and Urbanski
[13] and Lin \cite{lin} proved $f(K_3;4) = 15$. However, both upper bounds of Folkman
and of  Ne\v{s}et\v{r}il and R\"{o}dl for $f(K_3;3)$ are extremely large. Frankl and R\"odl \cite{frankl-rodl} first
gave a reasonable bound $f(K_3;3) \le 7 \times 10^{11}$.
Erd\H{o}s set a prize of \$100 for the challenge $f(K_3;3)\le 10^{10}$. This reward was claimed
by Spencer [10, 11], who proved that $f(K_3;3) < 3 \times 10^9$.
Erd\H{o}s then offered another \$100 prize (see \cite{chung-graham}, page 46) for the new challenge $f(K_3;3) < 10^6$.
Chung and Graham \cite{chung-graham} conjectured further $f(K_3;3)<10000$, which was confirmed by Lu \cite{lu} with  $f(K_3;3) < 9697$, and
by Dudek and  R\"{o}dl \cite{dud-r} with more computer aid.

Let us call a family ${\cal G}$ of graphs $G_n$ of order $n$ to be Ramsey linear if there exists a constant $c=c({\cal G})>0$ such that $R(G_n)\le cn$ for any $G_n\in {\cal G}$. Similarly, we call ${\cal G}$ to be Folkman $p$-linear if $f(G_n;p)\le cn$ for any $G_n\in {\cal G}$, where $p$ is a constant.  Let  $\Delta(G_n)$ be the maximum degree of $G_n$ of order $n$ and set a family of graphs as
\[
{\cal G}_{\Delta}=\{G_n\;|\;\Delta(G_n)\le \Delta\}.
\]
A result of Chv\'atal, R\"odl, Szemer\'edi and Trotter \cite{ch-ro-sz-tr} is as follows.
\begin{theorem}\label{r-linear}
The family ${\cal G}_{\Delta}$ is Ramsey linear.
\end{theorem}
The proof of Theorem \ref{r-linear} is a remarkable application of Szemer\'edi regularity lemma, in which they used the general form of the lemma.
In order to generalize Theorem \ref{r-linear} to Folkman number, we shall have a multi-partite regularity lemma as follows.
\begin{theorem}\label{new-2}
For any $\epsilon>0$ and integers $m\ge1$ and $p\ge 2$, there exists an $M =M(\epsilon,m,p)$
such that each $p$-partite graph $G(V^{(1)},\dots,V^{(p)})$ with $|V^{(s)}|\ge M$, $1\le s\le p$, has a partition
$\big\{V^{(s)}_1,\cdots, V^{(s)}_k\big\}$ for each $V^{(s)}$, where $k$ is same for each part $V^{(s)}$ and $m\le k \le M$, such that

(1) $\big||V^{(s)}_i|-|V^{(s)}_j|\big|\le 1$ for each $s$;

(2) All but at most $\epsilon k^2{p\choose2}$ pairs $\big(V^{(s)}_i,V^{(t)}_j\big)$, $1\le s<t\le p$, $1\le i,j \le k$, are $\epsilon$-regular.
\end{theorem}

Using the above Theorem \ref{new-2}, we can deduce the following result on the Folkman $p$-linearity of ${\cal G}_{\Delta}$ for some fixed $p$.

\begin{theorem}\label{f-linear}
Let $\Delta\ge 3$ be an integer and $p=R(K_{\Delta})$. Then the family ${\cal G}_{\Delta}$ is Folkman $p$-linear.
\end{theorem}

Note that for sub-family ${\cal G}_{\Delta,\chi}$ consisting of $G\in{\cal G}_{\Delta}$ with $\chi(G) \le \chi$, we can take $p=R(K_\chi)$ such that ${\cal G}_{\Delta,\chi}$ is Folkman $p$-linear. A natural problem is asking what is a smaller $p$ such that ${\cal G}_{\Delta}$ is Folkman $p$-linear.

For an integer $r\ge 2$, we call an edge coloring of a graph by $r$ colors as an $r$-edge coloring of the graph. For graphs $F$ and $G$, let $F\to (G)_r$ signify that any $r$-edge coloring of $F$ contains a monochromatic $G$. Thus $R_r(G)$ is the smallest $N$ such that $K_N\to (G)_r$, and $f_r(G;p)$ is the smallest $N$ such that there exists a graph $F$ of order $N$ with $\omega(F)=p$ satisfying $F\to (G)_r$. Theorem \ref{f-linear} can be generalized as follows.

\begin{theorem}\label{r-f-linear}
Let $\Delta\ge 3$ and $r\ge 2$ be integers and $p=R_r(K_{\Delta})$. Then, there is some constant $c=c(\Delta,r)>0$ such that $f_r(G_n,p)\le cn$ for any $G_n\in {\cal G}_{\Delta}$.
\end{theorem}

\section{Multi-partite regularity lemma}

Let $A$ be a set of positive integers and $A_n=A\cap \{1,\dots,n\}$. In the 1930s, Erd\H{o}s and Tur\'an conjectured that if  $\overline{\lim}\frac{|A_n|}{n} >0$, then $A$ contains  arbitrarily long arithmetic progressions. The conjecture in case of length $3$ was proved by Roth \cite{roth53,roth54}.
The full conjecture was proved by Szemer\'edi \cite{sze75} with a deep and complicated combinatorial argument. In the proof he used a result, which is now called the bipartite regularity lemma, and then he proved the general regularity lemma in \cite{sze78}. The lemma has become a totally new tool in extremal graph theory. Sometimes the regularity lemma is called uniformity lemma, see e.g., Bollob\'as \cite{bol98} and Gowers \cite{gow}. For many applications, we refer the readers to the survey of Koml\'os and Simonovits \cite{koml-simo}. In this note, we shall discuss multi-partite regularity lemma in slightly different forms.

Let $G(U,V)$ be a bipartite graph on two color classes $U$ and $V$. For $X\subseteq U$ and $Y\subseteq V$,
denote by  $e(X,Y)$  the number of edges between $X$ and $Y$ of $G$. The ratio \label{'edge-b-d'}
\[
d(X,Y) = \frac{e(X,Y)}{|X||Y|}
\]
is called the edge density of $(X,Y)$, which is the probability that any pair $(x,y)$ selected randomly from $X\times Y$ is an edge.
Clearly $0\le d(X,Y)\le 1$.

The first form of regularity lemma given by Szemer\'edi in \cite{sze75} is as follows, in which corresponding to each subset $U_i$ in the partition of $U$, we have to choose its own partition $V_{i,j}$ of $V$.

\begin{lemma}[Bipartite Regularity Lemma-Old Form]\label{old}
For any positive $\epsilon_1,\epsilon_2,\delta, \rho_1,\rho_2$, there exist $k_1, k_2$, $M_1, M_2$ such that every bipartite graph $G(U,V)$ with $|U|>M_1$ and $|V|>M_2$, there exist disjoint $U_i\subset U$, $i<k_1$, and for each $i<k_1$, disjoint $V_{i,j}\subset V$, $j<k_2$, such that:

(1) $|U-\cup_{i<k_1}U_i|<\rho_1 |U|$, and $|V-\cup_{j<k_2}V_{i,j}|<\rho_2 |V|$ for any $i<k_1$;

(2) For all $i<k_1$, $j<k_2$, $X\subseteq U_i$ and $Y\subseteq V_{ij}$ with $|X|>\epsilon_1|U_i|$ and $|Y|>\epsilon_2|V_j|$, we have
\[
d(X,Y)\ge d(U_i,V_{ij})-\delta;
\]

(3) For all $i<k_1$, $j<k_2$ and $x\in U_i$, $|N(x)\cap V_{i,j}|\le (d(U_i,V_{ij})+\delta)|V_{i,j}|$.
\end{lemma}

For $\epsilon>0$, a disjoint pair $(X,Y)$ is called $\epsilon$-{\em regular} if any $X'\subseteq X$ and
$Y'\subseteq Y$ with $|X'| > \epsilon |X|$ and $|Y'| > \epsilon
|Y|$ satisfy
\[
| d(X,Y) - d(X',Y') | \le \epsilon.
\]

We shall call $U_0=U-\cup_{i<k_1}U_i$, and $V_0=V-\cup_{j<k_2}V_{i,j}$ in Theorem \ref{old} to be the {\em exceptional sets}. The following is the general regularity lemma of Szemer\'edi \cite{sze78}, in which the partition $C_0,C_1,\dots,C_k$ is {\em equitable} in sense of that all sets $C_i$ other than the exceptional set $C_0$ have the same size.

\begin{lemma}[General Regularity Lemma]\label{general}
For any $\epsilon>0$ and any $m\ge 1$, there exists $M =M(\epsilon,m)>m$
such that every graph $G$ of order at least $m$ has a partition $C_0, C_1,\dots,C_k$  with  $m\le k \le M$ such that

(1) $|C_1|=|C_2|=\cdots=|C_k|$ and $|C_0|\le \epsilon n$;

(2) All but at most $\epsilon k^2$ pairs $(C_i,C_j)$ with $1\le i< j \le k$ are $\epsilon$-regular.
\end{lemma}

There are many generalizations of Szemer\'edi regularity lemma, in particular, Frankl and R\"{o}dl \cite{frankl-rodl-1992} generalized it to hypergraphs and later Chung \cite{chung} formulated regularity lemma on $t$-uniform hypergraphs when discussing the problems of quasi-random hypergraphs.

The regularity lemma has numerous applications in various areas, mainly in extremal graph theory such as \cite{ch-ro-sz-tr} by Chv\'atal, R\"odl, Szemer\'edi and Trotter. In an application, Eaton and R\"{o}dl \cite{eaton-rod} obtained a form of the regularity lemma for $p$-partite $p$-uniform hypergraph. To state their result for multi-partite graph, let us have some definitions.

Let $G(V^{(1)},\dots,V^{(p)})$ be a $p$-partite graph on vertex set $\cup_{i=1}^p V^{(i)}$. Consider partitions of the set
$V^{(1)}\times \cdots\times V^{(p)}$, where each partition class is of the form $W_1\times \cdots\times W_p$, $W_i \subseteq V^{(i)}$,
$1\le i \le p$, which is called {\em cylinders}. Let us say that a cylinder $W_1\times\cdots\times W_p$ is
$\epsilon$-regular if the subgraph of $G$ induced on the set $\cup_{i=1}^p W_i$ is such that all pairs $(W_i,W_j)$,
$1\le i < j \le p$, are $\epsilon$-regular.

Eaton and R\"{o}dl stated their result with exceptional $p$-tuples instead of exceptional sets, for which Alon, Duke, Lefmann, R\"odl and Yusterk studied the computational difficulty of finding such a regular partition in \cite{a-d-l-r-y}.

\begin{lemma} \label{eaton-r}
Let $G(V^{(1)},\dots,V^{(p)})$ be a $p$-partite graph with $|V^{(i)}| = n$, $i = 1,\dots,p$.
Then for every $\epsilon > 0$ there exists a partition of $V^{(1)}\times\cdots\times V^{(p)}$ into $k$ cylinders with $k \le 4^h$, where $h =\frac{{p\choose 2}}{\epsilon^5}$, such that all but at most $\epsilon n^p$ of the $p$-tuples $(v_1,\dots,v_p)$ of $V^{(1)}\times\cdots\times V^{(p)}$ are in $\epsilon$-regular cylinders of the partition.
\end{lemma}

Note that in Lemma \ref{eaton-r}, the transverse  section $\{W_i\}$ of the partition is a partition of $V^{(i)}$, which may be not equitable, and for $i\neq j$, the numbers of subsets in the partitions $\{W_i\}$  and $\{W_j\}$  may be different. We shall have a multi-partite regularity lemma as follows.

\begin{lemma}\label{new-1}
For any $\epsilon>0$ and integers $m\ge1$ and $p\ge2$, there exists $M =M(\epsilon,m,p)$
such that each $p$-partite graph $G(V^{(1)},\dots,V^{(p)})$ with $|V^{(s)}|\ge M$, $1\le s\le p$, has a partition
$\big\{V^{(s)}_0, V^{(s)}_1,\cdots, V^{(s)}_k\big\}$ for each $V^{(s)}$, where $k$ is same for each part $V^{(s)}$ and $m\le k \le M$, such that

(1) $\big|V^{(s)}_1\big|=\cdots=\big|V^{(s)}_k\big|$  and  $\big|V^{(s)}_0\big|\le \epsilon \big|V^{(s)}\big|$ for each $s$;

(2) All but at most $\epsilon k^2{p\choose2}$ pairs $\big(V^{(s)}_i,V^{(t)}_j\big)$, $1\le s< t\le p$, $1\le i,j \le k$, are $\epsilon$-regular.
\end{lemma}

The following multicolor multi-partite regularity lemma is an analogy of Theorem \ref{new-2}, which is needed for proof of Theorem \ref{r-f-linear}.

\begin{lemma}\label{multi-color-reg-le}
For any $\epsilon>0$ and integers $m\ge1$, $p\ge 2$ and $r\ge1$, there exists an $M =M(\epsilon,m,p,r)$
such that if the edges of a $p$-partite graph $G(V^{(1)},\dots,V^{(p)})$ with $|V^{(s)}|\ge M$, $1\le s\le p$ are $r$-colored, then all monochromatic graphs have the same partition $\big\{V^{(s)}_1,\cdots, V^{(s)}_k\big\}$ for each $V^{(s)}$, where $k$ is same for each part $V^{(s)}$ and $m\le k \le M$, such that

(1) $\big||V^{(s)}_i|-|V^{(s)}_j|\big|\le 1$ for each $s$;

(2) All but at most $\epsilon k^2 r{p\choose2}$ pairs $\big(V^{(s)}_i,V^{(t)}_j\big)$, $1\le s<t\le p$, $1\le i,j \le k$, are $\epsilon$-regular in each monochromatic graph.
\end{lemma}

\section{Proofs for multi-partite regularity lemma}

In this section, we prove Lemma \ref{new-1},  Theorem \ref{new-2} and Lemma \ref{multi-color-reg-le}. To reduce the complicity of notations in the proofs, we shall prove them in case $p=2$, which are bipartite regularity lemmas.

\begin{lemma}\label{le-cont}
Let $G(U,V)$ be a bipartite graph and let $X\subseteq U$ and $Y\subseteq V$. If
$X'\subseteq X$ and $Y'\subseteq Y$ satisfy $|X'| > (1-\delta)|X|$
and $|Y'|  > (1-\delta)|Y|$, then
\[
|d(X',Y')-d(X,Y)| < 2\delta \hspace{3mm} \mbox{and} \hspace{3mm} |d^2(X',Y')-d^2(X,Y)| < 4\delta.
\]
\end{lemma}

A crucial point for the regularity lemma is that the number $k$ of
classes in partition is bounded for any graph. For proofs, we need the well-known defect form of Cauchy-Schwarz inequality.

\begin{lemma}\label{cau-sch}
Let  $d_i$ be reals and $s> t \ge 1$ be integers. If
\[
\frac{1}{s} \sum_{i=1}^s d_i =\frac{1}{t} \sum_{i=1}^t d_i + \delta,
\]
then
\[
\frac{1}{s} \sum_{i=1}^s d_{i}^2 \ge
\left( \frac{1}{s} \sum_{i=1}^s d_i \right)^2  +\frac{t\delta^2}{s-t}.
\]
\end{lemma}

Let $G(U,V)$ be a bipartite graph, a partition
\[
{\cal P}=\Big\{U_i,V_j\Big| \;0\le i,j\le k\Big\},
\]
where $U=\cup_{p=1}^k U_i$ and $V=\cup_{p=1}^k V_i$, is called to be an {\em equitable} partition of $U\cup V$ with exceptional classes $U_0$ and $V_0$ if $|U_i|=|U_j|$ and $|V_i|=|V_j|$ for $1\le i, j\le k$. For convenience, we say an equitable partition ${\cal P}$ is {\em $\epsilon$-regular} if all but at most $\epsilon k^2$ pairs of $(U_i,V_j)$ are $\epsilon$-regular.
Define
\[
q({\cal P}) = \frac{1}{k^2} \sum_{1\le i,j \le k} d^2(U_i,V_j).
\]
It is easy to see that $0\le q({\cal P}) \le 1$ since $0\le d(U_i,V_j) \le 1$.

In the following, we will show that if ${\cal P}$ is not
$\epsilon$-regular, then there is a partition ${\cal P}'$ with the new
exceptional classes a bit larger than the old one, but $q({\cal P}')
\ge q({\cal P}) +\frac{\epsilon^5}{4}$. Do this again if ${\cal P}'$ is
not $\epsilon$-regular yet. The number of iterations is thus at
most $4/\epsilon^5$ in order to obtain an
$\epsilon$-regular partition. Without loss of generality, we
assume that $0<\epsilon \le 1/2$ since if $\epsilon>1/2$, one can
take $M(\epsilon,m)$ to be $M(1/2,m)$.

\begin{lemma}\label{main-le}
Let $G(U,V)$ be a bipartite graph with $|U|=n_1\ge M$ and $|V|=n_2\ge M$, which has an equitable partition
\[
{\cal P}=\Big\{U_i,V_j\Big| \;0\le i,j\le k\Big\}
\]
with exceptional classes  $U_0$ and $V_0$. Suppose $2^k \ge 16/\epsilon^5$,
$
|U_i|=c_1 \ge 2^{3k}\; \mbox{and} \;\; |V_j|=c_2 \ge 2^{3k}.
$
We have if ${\cal P}$ is not $\epsilon$-regular, then there is an equitable partition
\[
{\cal P}'=\Big\{U'_i,V'_j \Big| \,0\le i,j\le \ell\Big\}
\]
with exceptional class $U'_0\supseteq U_0$ and $V'_0\supseteq V_0$, and  $\ell = k(4^k-2^k)$ satisfying

(1) \ $|U'_0| \le |U_0| + n_1/2^{k-1}$ and $|V'_0| \le |V_0| + n_2/2^{k-1}$;

(2) \ $q({\cal P}') \ge q({\cal P}) + \epsilon^5/4$.
\end{lemma}
{\bf Proof. } Separate all pairs $(i,j)$, $1\le i,j\le k$,
of indices into $S$ and $T$, corresponding with that the pair
$(U_i,V_j)$ is $\epsilon$-regular or not, respectively. For
$(i,j)\in S$, set $U_{ij} =V_{ji} =\emptyset$, and for  $(i,j)\in
T$, set $U_{ij} \subseteq U_i$ and $V_{ji} \subseteq V_j$ with
$|U_{ij}|>\epsilon c_1$, $|V_{ji}| > \epsilon c_2$, and
\[
|d(U_{ij},V_{ji}) -d(U_i,V_j)| > \epsilon.
\]

For fixed $i$, $1\le i \le k$, consider an equivalence relation
$\equiv$ on $U_i$ as $x\equiv y$ if and only if both $x$ and $y$
belong to the same ${U_{ij}}'$s. The equivalence classes are
atoms of algebra induced by $U_{ij}$, and each $U_i$
has at most $2^{k}$ atoms. Similarly, each $V_j$
has at most $2^{k}$ atoms.

For $p=1,2$, set $d_p=\lfloor c_p/4^k \rfloor$.
Let us cut each atom in $U_i$ into pairwise disjoint $d_1$-subsets.
Denote by $z$ for the maximal number of these $d_1$-subsets that one can take, clearly $z\ge 4^k - 2^{k}$ as $zd_1+2^{k}(d_1-1) \ge c_1$. Set
\[
H=4^k - 2^{k},
\]
and take
{\em exactly} $H$ such $d_1$-subsets and add the remainder to the
``rubbish bin" to get a new exceptional set $U'_0$. Label all these
$d_1$-subsets in $U_i$ as $D_{i1},\dots,D_{iH}$. Set $U'_0=U_0
\cup \Big[ \cup_{i=1}^k \Big( U_i \setminus \cup_{h=1}^H D_{ih}
\Big)\Big]$, and so $|U'_0| =|U_0| +k(c_1-Hd_1)$. Since
\[
Hd_1 \ge (4^k -2^{k}) (\frac{c_1}{4^k} -1) > c_1 - \frac{c_1}{2^{k-1}}
\]
by noting $c_1\ge 2^{3k}$, we have $|U'_0| \le
|U_0|+n_1/2^{k-1}$. Rename $D_{ih}$ as $U'_s$ for $1\le s
\le \ell$, where $\ell =kH$.

Similarly, we can cut each atom in $V_j$ into pairwise disjoint $d_2$-subsets and take $H$ such subsets $E_{j1},\dots,E_{jH}$ in $V_j$.
Set $V'_0=V_0\cup \Big[ \cup_{i=1}^k \Big( V_j \setminus \cup_{h=1}^H E_{jh}\Big)\Big]$, and similarly $|V'_0| \le|V_0|+n_2/2^{k-1}$. Rename $E_{jh}$ as $V'_t$ for $1\le t\le \ell$.

Denote the new equitable partition by
\[
{\cal P}'=\Big\{U'_i,V'_j \Big| \,0\le i,j\le \ell\Big\}
\]
of $U\cup V$ with exceptional classes $U'_0\supseteq U_0$ and $V'_0\supseteq V_0$. All that remains is to show $q({\cal P}') \ge q({\cal P}) +\epsilon^5/4$.

For $1\le i,j\le k$, set
\begin{align*}
\overline{U_i} &= \cup_{h=1}^H D_{ih}, \hspace{3mm}
                      \overline{U_{ij}} = \cup \{ D_{ih}: D_{ih} \subseteq U_{ij} \}, \hspace{3mm}  \mbox{and} \hspace{3mm}
\\
\overline{V_j} &= \cup_{h=1}^H E_{jh}, \hspace{3mm}
                      \overline{V_{ji}} = \cup \{ E_{jh}: D_{jh} \subseteq V_{ji} \}.
\end{align*}
Set partition $\overline{{\cal P}}=\{ U'_0, \overline{U_1},\dots,\overline{U_k};  V'_0, \overline{V_1},\dots,\overline{V_k}\}$
with exceptional class $U'_0$ and $V'_0$.

{\em Claim 1.}  $\; q(\overline{{\cal P}}) \ge q({\cal P}) -
\epsilon^5/2$.

{\em Proof of Claim 1.}  \ Note that $\frac{|U_i\setminus \overline{U_i}|}{|U_i|}<\frac{1}{2^{k-1}}<\frac{\epsilon^5}{8}$ and $\frac{|V_j\setminus \overline{V_j}|}{|V_j|}<\frac{\epsilon^5}{8}$ for any pair $(U_i,V_j)$,
we have
\begin{equation}  \label{eq-2}
|d(\overline{U_i},\overline{V_j}) - d(U_i,V_j) | \le \frac{\epsilon^5}{4}
\end{equation}
 by Lemma \ref{le-cont}.
Hence $d^2(\overline{U_i},\overline{V_j}) \ge d^2(U_i,V_j) - \epsilon^5/2$, which implies that
$q(\overline{{\cal P}}) \ge q({\cal P}) - \epsilon^5/2$ as claimed.

\smallskip

{\em Claim 2.}  \ If  $(i,j)\in T$, then
$|d(\overline{U_{ij}},\overline{V_{ji}})-d(\overline{U_i},\overline{V_j})|
> \frac{15}{16}\epsilon$.

{\em Proof of Claim 2.} \ Clearly,  $\frac{|U_{ij}\setminus \overline{U_{ij}}|}{|U_{ij}|} \le
  \frac{|U_i\setminus \overline{U_i}|}{|U_i|}\frac{|U_i|}{|U_{ij}|} \le \frac{\epsilon^4}{8}$
and $\frac{|V_{ji}\setminus \overline{V_{ji}}|}{|V_{ji}|}\le \frac{\epsilon^4}{8}$,
which and Lemma \ref{le-cont} give
\begin{align}\label{eq-4}
|d(\overline{U_{ij}},\overline{V_{ji}}) - d(U_{ij},V_{ji}) | \le \frac{\epsilon^4}{4}.
\end{align}
Therefore, if $(i,j)\in T$, the bounds (\ref{eq-2}) and (\ref{eq-4})
with the fact that $0<\epsilon\le 1/2$ will yield the desired inequality.

\vspace{1mm}

Let us return to the partition ${\cal P}'$ in which each class is either a
$d_1$-subset $D_{iu}$ or a $d_2$-subset $E_{jv}$ except $U'_0$ and $V'_0$. For any pair $(U_i,V_j)$,
\[
d(\overline{U_i},\overline{V_j}) = \frac{1}{H^2}\sum_{1\le u, v\le H} d(D_{iu},E_{jv})
\]
since $|\overline{U_i}|=Hd_1$ and $|\overline{V_j}| =Hd_2$. Set
\[
A(i,j) = \frac{1}{H^2}\sum_{1\le u, v\le H} d^2 (D_{iu},E_{jv}).
\]
Then from Cauchy-Schwarz inequality, for any pair $(i,j)$, we have
\begin{align}\label{ijs}
A(i,j) \ge d^2 (\overline{U_i},\overline{V_j}).
\end{align}

If $(i,j)\in T$, we have some gain. Let $R=R(i,j)$ be the set of
indices $(u,v)$ such that $D_{iu}\in \overline{U_{ij}}$ and
$E_{jv}\in \overline{V_{ji}}$. Then
\[
d(\overline{U_{ij}},\overline{V_{ji}}) = \frac{1}{|R|} \sum_{(u,v)\in R} d(D_{iu},E_{jv}).
\]
Note that
$\frac{|R|}{H^2}=\frac{{\overline{U_{ij}}}\,\,{\overline{V_{ji}}}}{{\overline{U_{i}}\,\,\overline{V_{j}}}}\ge\big((1-2^{-7})\epsilon\big)^2,
$
So Lemma \ref{cau-sch} and Claim 2 imply
\begin{align}\label{ijt}
A(i,j)  \ge  d^2 (\overline{U_i},\overline{V_j}) +
 \frac{|R|}{H^2} \left(d (\overline{U_i},\overline{V_j})-d(\overline{U_{ij}},\overline{V_{ji}}\right)^2
  \ge  d^2 (\overline{U_i},\overline{V_j}) + \frac{3}{4}\epsilon^4.
\end{align}
Noticing that $\ell =kH$, we have
\begin{align*}
q({\cal P}') &= \frac{1}{\ell^2} \sum_{1\le s,t \le \ell} d^2 (U'_s,V'_t) = \frac{1}{k^2} \frac{1}{H^2}\sum_{1\le i, j\le k} \sum_{1\le u, v\le H} d^2 (D_{iu},E_{jv})= \frac{1}{k^2} \sum_{1\le i, j\le k} A(i,j).
\end{align*}
Now, combine inequalities (\ref{ijs}) and (\ref{ijt}), and recall Claim 1 and that ${\cal P}$ is not $\epsilon$-regular, we have
\begin{align*}
q({\cal P}')   & \ge   \frac{1}{k^2} \left[ \sum_{(i,j)\in S} d^2 (\overline{U_i},\overline{V_j})
  +  \sum_{(i,j)\in T}\left( d^2 (\overline{U_i},\overline{V_j}) + \frac{3}{4}\epsilon^4 \right)\right]
                \ge q({\cal P}) + \frac{\epsilon^5}{4}.
\end{align*}
This completes the proof of Lemma \ref{main-le}.            \hfill    $\Box$

\medskip

\noindent{\bf Proof of Lemma \ref{new-1}.}
Let $k_0$ be an integer such that $k_0\ge m$ and $2^{-k_0} \le
\epsilon^5/16$, and define $k_{i+1} = k_i(4^{k_i}-2^{k_i})$.  Set
$M_i=k_i 2^{3k_i}$ and $M=M_t$.  Lemma \ref{main-le} implies that at most $t=4\lfloor
\epsilon^{-5}\rfloor$ iterations will yield a required partition, which completes the proof of Lemma \ref{new-1}.
 \hfill    $\Box$

\medskip
\noindent
{\bf Proof of Theorem \ref{new-2}.}  \ For given $\epsilon>0$ and $m\ge 1$, Theorem \ref{new-2} implies that there is an $M>m$
and an equitable and $\frac{\epsilon^2}{4}-$regular partition ${\cal P}=\{U_i,V_j|0\le i,j\le k\}$ with $m\le k \le M$. Since $|U_0|<\frac{\epsilon^2}{4}n_1$, we have $\lfloor (1-\epsilon^2/4)n_1/k \rfloor \le |U_i| \le n_1/k$. Partition $U_0$ into
$k$ classes $U_{01},U_{02},\dots,U_{0k}$ such that $|U_{0i}|=\lfloor |U_0|/k\rfloor$ or $|U_{0i}|=\lceil |U_0|/k \rceil$. Set $U'_i=U_i \cup U_{0i}$, clearly $|U'_i|=\lfloor n_1/k \rfloor$ or $|U'_i|=\lceil n_1/k \rceil$.
Similarly, let us partition $V_0$ into
$k$ classes $V_{01},V_{02},\dots,V_{0k}$ such that $|V_{0i}|=\lfloor |V_0|/k\rfloor$ or $|V_{0i}|=\lceil |V_0|/k \rceil$. Set $V'_i=V_i \cup V_{0i}$, we have the sizes of any $V'_i$ and $V'_j$ differ at most by one. Then the Partition ${\cal P}'=\{U'_i,V'_j|0\le i,j\le k\}$ is as desired by noting that if a pair $(U_i,V_j)$ is $\frac{\epsilon^2}{4}-$regular, then $(U'_i,V'_j)$ is $\epsilon$-regular. \hfill  $\Box$

\medskip

\noindent{\bf Proof of Lemma \ref{multi-color-reg-le}.}
A similar proof as Theorem \ref{new-2}, but modify the definition of index by summing the indices for each color,

\[
q({\cal P}) = \frac{1}{k^2} \sum_{1\le h \le r}\sum_{1\le s<t \le p}\sum_{1\le i,j \le k} d^2\big(V^{(s)}_i,V^{(t)}_j\big).
\]
Then we have analogy of Lemma \ref{new-1} for multi-color case. Furthermore, we have Lemma \ref{multi-color-reg-le}. \hfill $\Box$

\section{A Folkman linear family}

In this section, we shall apply multi-partite regularity lemma to the Folkman numbers involving the family ${\cal G}_{\Delta}$ of graphs with maximum degree bounded. In order to prove Theorem \ref{f-linear} and Theorem \ref{r-f-linear}, we shall establish the following Lemma, in which $K_p(k)$ is the complete $p$-partite graph with $k$ vertices in each part.

\begin{lemma}\label{turan}
For integers $k\ge 1$ and $p\ge 2$, let $t_p(k)$ be the maximum number of edges in a subgraph of $K_p(k)$ that contains no $K_p$. Then
\[
t_p(k)= \left[{p\choose 2}-1\right] k^2.
\]
\end{lemma}
{\bf Proof.} By deleting all edges between a pair of parts of $K_p(k)$, we have the lower bound for $t_p(k)$ as required. On the other hand, we shall prove  by induction of $k$ that if a subgraph $G=G(V^{(1)},\dots,V^{(p)})$ of $K_p(k)$ contains no $K_p$, then $e(G)\le \left[{p\choose 2}-1\right] k^2$. Suppose $k\ge 2$ and $p\ge3$ as it is trivial for $k=1$ or $p=2$. Furthermore, suppose  that
$G$ has the maximum possible number of edges subject to this condition. Then $G$ must contain $K_{p}-e$ as a
subgraph, otherwise we could add an edge and the resulting graph would still not contain $K_p$. Pick a vertex set $X$ consisting of a vertex from each $V^{(i)}$ for $i=1,2,\dots,p$ such that $e(X)$ is maximum among all such vertex subsets, and so $e(X)={p\choose 2}-1=\frac{(p+1)(p-2)}{2}$. We may suppose that $X$ induces a complete graph of order $p$ with an edge $v_1v_2$ missing, where $v_1\in V_1$ and $v_2\in V_2$. Let $Y=V(G)\setminus X$, clearly each part of $Y$ has $k-1$ vertices. Now, by noticing the fact that no vertex in $V^{(i)}\cap V(Y)$ is adjacent to all the vertices of $X\setminus\{v_i\}$ for $i=1,2$ since $G$ contains no $K_p$, we can safely deduce the desired upper bound of $t_p(k)$ by a simple calculation, which completes the induction hypothesis hence the proof.
  \hfill  $\Box$

\begin{lemma}\label{pre}
 Let $(A,B)$ be an $\epsilon$-regular pair of density $d\in
(0,1]$, and $Y\subseteq B$ with $|Y|\ge \epsilon |B|$. Then there
exists a subset $A'\subseteq A$ with $|A'|\ge (1-\epsilon)|A|$,
each vertex in $A'$ is adjacent to at least $(d-\epsilon)|Y|$
vertices in $Y$.
\end{lemma}
{\em Proof.} Let $X$ be the set of vertices with fewer than
$(d-\epsilon)|Y|$ neighbors in $Y$. Then $e(X,Y)
<(d-\epsilon)|X||Y|$, so $d(X,Y)<d-\epsilon$. Since $(A,B)$ is
$\epsilon$-regular, this implies that $|X|<\epsilon |A|$. \hfill
$\Box$

\medskip

\noindent{\bf Proof of Theorem \ref{f-linear}.}  We will consider a red/blue edge coloring of $K_p(cn)$. Denote by $H_R$ and $H_B$ the subgraphs spanned by red edges and blue edges, respectively. Note that a partition obtained by applying Theorem \ref{new-2} for $H_R$ is such a partition for $H_B$.

Let $p=R(K_{\Delta})$ as defined.  Clearly, we can only consider graphs $G=G_n$ in ${\cal G}_{\Delta}$ with $n\ge \Delta+2$. Choose $\epsilon=\min\{\frac{1}{p^2},\frac{1}{m}\}$, where $m$ is a positive integer such that
\[
(1-\Delta\epsilon)(1/2-\epsilon)^{\Delta}\,m \ge 1 \hspace{3mm}\mbox{hence}\hspace{3mm} (1-\Delta\epsilon)(1/2-\epsilon)^{\Delta} \ge\epsilon.
\]
Let $M=M(\epsilon,m, p)>2m$ be the integer determined by $\epsilon$
and $p$ in  Theorem \ref{new-2} for $H_R$. Finally, let $c=mM$ which is a
constant determined completely by $\Delta$. We shall show that either $H_R$ contains $G$ or $H_B$ contains $G$, hence $f(G;p)\le cpn$.

Let the vertex set of the $K_{p}(cn)$ be $V=V^{(1)}\cup\cdots\cup V^{(p)}$ with $|V_\ell|=cn$ for $1\le \ell\le p$. There is
a partition of $V$, in which each $V^{(\ell)}$ is partitioned into $\big\{V^{(\ell)}_1,\dots,V^{(\ell)}_k\big\}$  with
$\big||V^{(\ell)}_i|-|V^{(\ell)}_j|\big|\le 1$ and $m\le k \le M$, and all but at
most $\epsilon k^2{p\choose 2}$ pairs $\big(V^{(s)}_i,V^{(t)}_j\big)$, $1\le i,j\le k$, $1\le s\neq t\le p$, are
$\epsilon$-regular.

Let $F$ be the subgraph of $K_p(k)$, whose vertices are $\big\{V^{(\ell)}_i\;|\;1\le \ell\le p, 1\le i\le k\big\}$ in which a pair  $\big(V^{(s)}_i,V^{(t)}_j\big)$ for $s\neq t$ is adjacent if and only if
the pair  is $\epsilon$-regular in $H_R$.  Then the number of edges of
$F$ is at least
\[
(1-\epsilon)k^2{p\choose 2} >\left[{p\choose 2}-1\right]k^2=t_p(k).
\]
By Lemma \ref{turan}, $F$ contains a complete graph $K_p$. Without loss of generality, assume that $V^{(1)}_1,\dots, V^{(p)}_1$ are pairwise $\epsilon$-regular.
Color an edge between a pair $(V^{(s)}_1, V^{(t)}_1)$ green  if $d\big(V^{(s)}_1,V^{(t)}_1\big)\ge 1/2$, or white if $d\big(V^{(s)}_1,V^{(t)}_1\big)< 1/2$. As $p=R(K_{\Delta})$, we have $\Delta$ sets in $\{V^{(1)}_1,V^{(2)}_1, \dots, V^{(p)}_1\}$ such that they form a monochromatic $K_{\Delta}$.
We may assume that the color is green since otherwise we consider the graph $H_B$.

Relabeling the sets in the partition if necessary, we assume that
$V^{(1)}_1,V^{(2)}_1, \dots, V^{(\Delta)}_1$ are pairwise $\epsilon$-regular in $H_R$, and $d\big(V^{(s)}_1,V^{(t)}_1\big) \ge 1/2$.
Write
\[
C_1=V^{(1)}_1,C_2=V^{(2)}_1, \dots, C_{\Delta}=V^{(\Delta)}_1.
\]
Note that if $Y_i \subseteq C_i$ with $|Y_i| \ge
(1-\Delta\epsilon)(1/2-\epsilon)^{\Delta}|C_i|$,
 then $|Y_i| \ge\epsilon |C_i|$,
which is the preparation for using Lemma \ref{pre}, and
\[
 |Y_i|\ge (1-\Delta\epsilon)(1/2-\epsilon)^{\Delta} \frac{cn}{M} \ge n,
 \]
which will give us enough room to maneuver for constructing a color class of $G$.

Note that if a graph is neither a complete graph nor an odd cycle, then its chromatic number is at most $\Delta(G)$. For considered graph $G=G_n$, as $n\ge \Delta+2$ and $\Delta\ge 3$, we have $\chi(G)\le \Delta$.

Assume that $V(G)=\{u_1,u_2,\dots,u_n\}$. We shall show that the red
graph $H_R$ contains $G$ as a subgraph. We will choose
$v_1,v_2,\dots,v_n$ from the sets $C_1,\dots,C_{\Delta}$. Since  $\chi(G)\le \Delta$, so $V(G)$ can be
partitioned into $\Delta$ color classes, which defines a map
$\phi$: $\{1,\dots,n\}\to \{1,\dots,\Delta\}$, where $\phi(i)$ is
the color of vertex $u_i$. Our aim is to define an embedding $u_i\to
v_i\in C_{\phi(i)}$, such that $v_iv_j$ is an edge of $H_R$ whenever
$u_iu_j$ is an edge of $G$.

Our plan is to choose the vertices $v_1,\dots,v_n$ inductively.
Throughout the induction, we shall have a target set $Y_i\subseteq
C_{\phi(i)}$ assigned to each $i$. Initially, $Y_i$ is the entire
$C_{\phi(i)}$. As the embedding proceeds, $Y_i$ will get smaller and
smaller. Some vertices will be deleted in procedure. But any
$C_{\phi(i)}$ will really have some vertices deleted at most
$\Delta$ times. To make this approach work, we have to ensure $Y_i$
do not get too small.

Let us begin the initial step. Set
\[
Y_1^0=C_{\phi(1)},\; Y_2^0=C_{\phi(2)}, \; \dots,\; Y_n^0
=C_{\phi(n)}.
\]
Note that $Y_i^0$ and $Y_j^0$ are not necessarily distinct sets.

We then begin the first step by considering $u_1$, for which $v_1$
will be selected from $Y_1^0$, and its neighbors,
$u_\alpha,\dots,u_\beta$, say. Suppose that the degree of $u_1$ is
$d$. By using Lemma \ref{pre} repeatedly, we know that there exists
a subset $Y_1^1\subseteq Y_1^0$ with $|Y_1^1|\ge
(1-d\epsilon)|Y_1^0|\ge n$, such that each vertex in $Y_1^1$ has at
least $(1/2-\epsilon)|Y_{j}^0|$ neighbors in $Y_{j}^0$, where
$j=\alpha, \dots, \beta$. Choose an arbitrary vertex
 $v_1$ from $Y_1^1$. For $j=\alpha,\dots,\beta$,
define $Y_j^1$ be the neighborhood of $v_1$  in $Y_j^0$. For $j\ge
2, j\not=\alpha,\dots,\beta$, define $Y_j^1=Y_j^0$, that is, no
vertices are deleted from such $Y_j^0$. In this step, $v_1$ has
been chosen and it completely adjacent to $Y_j^1$ in $H$ whenever
$u_1$ and $u_j$ are adjacent in $G$.

In a general step, we consider $u_i$ and its neighbors. We will
choose $v_i$ for $u_i$ from $Y_{i}^{i-1}$. Suppose that $u_i$ has
$d_1$ neighbors in $\{u_1,\dots,u_{i-1}\}$, and $d_2$ neighbors,
$u_\alpha,\dots,u_\beta$, say, in $\{u_{i+1},\dots,u_n\}$. Then
$d_1+d_2\le \Delta$, and $|Y_i^{i-1}|\ge
(1/2-\epsilon)^{d_1}|Y_i^0|$. That is to say, the current set
$Y_i^{i-1}$ are obtained from $Y_i^0$ by deleting some vertices
$d_1$ times before this step. By using Lemma \ref{pre} repeatedly
again, we have a subset $Y_i^i \subseteq Y_i^{i-1}$ with $|Y_i^i|\ge
(1-d_2\epsilon)|Y_i^{i-1}|$ so that each vertex in $Y_i^i$ has at
least $(1/2-\epsilon)|Y_j^{i-1}|$ neighbors in $Y_j^{i-1}$, where
$j=\alpha,\dots,\beta$. Since
\begin{eqnarray*}
&& |Y_i^i| \ge (1-d_2\epsilon)|Y_i^{i-1}|  \ge
(1-d_2\epsilon)(1/2-\epsilon)^{d_1}|Y_i^{0}| \\
& & \ge  (1-\Delta\epsilon)(1/2-\epsilon)^{\Delta} |C_i|  \ge  n,
\end{eqnarray*}
we can choose a vertex $v_i$ from $Y_i^i$, which is distinct from
$v_1,\dots,v_{i-1}$ that have been chosen before this step, and some
may be from $Y_i^i$. For $j=\alpha,\dots,\beta$, define $Y_j^i$ to
be the neighborhood of $v_i$ in $Y_j^{i-1}$. For $j\ge i+1,
j\not=\alpha,\dots,\beta$, define $Y_j^i=Y_j^{i-1}$, that is, no
vertices are deleted from such $Y_j^{i-1}$.   Note that $v_i$ is
adjacent to any $v_j$, where $j<i$ and $u_j$ is adjacent to $u_i$,
and $v_i$ is completely connected with each set $Y_j^i$, in which a
neighbor of $v_i$ will be selected after this step.

It is easy to check that the condition for using Lemma \ref{pre}
can be satisfied since
 $(1-\Delta\epsilon)(1/2-\epsilon)^{\Delta}\ge \epsilon$.
We thus finished the general step hence the proof of Theorem \ref{f-linear}. \hfill  $\Box$

\medskip

\noindent{\bf Proof of Theorem \ref{r-f-linear}.}  For $p=R_r(K_{\Delta})$, take
$\epsilon=\min\{\frac{1}{p^2},\frac{1}{m}\}$, where $m$ is an integer such that
\[
(1-\Delta\epsilon)(1/r-\epsilon)^{\Delta}\,m \ge 1.
\]
In the proof, we use Lemma \ref{multi-color-reg-le}. We shall have $p$ sets, say $V_1^{(1)},\dots, V_1^{(p)}$, such that every pair $(V_1^{(s)},V_1^{(t)})$, $1\le s<t\le p$, is $\epsilon$-regular in each monochromatic graph. Connecting this pair with color $\ell$ if its edge density is at least $1/r$ in the monochromatic graph in color $\ell$, $1\le\ell\le r$. Then we have a $r$-edge coloring of $K_p$, which implies a monochromatic $K_{\Delta}$ in some color, say the color $a$. Hence we obtain $\Delta$ sets, say $V_1^{(1)},\dots, V_1^{(\Delta)}$, such that each pair $(V_1^{(s)},V_1^{(t)})$, $1\le s<t\le \Delta$, is $\epsilon$-regular in monochromatic graph of color $a$, and the edge density of the pair is at least $1/r$ in this color. The remaining proof is similar to that for Theorem \ref{f-linear}.  \hfill $\Box$

\end{document}